\documentclass[11pt]{amsart}
\usepackage{graphicx}
\usepackage{amscd}
\usepackage{amsmath}
\usepackage{amsfonts}
\usepackage{amssymb}
\newtheorem{theorem}{Theorem}[section]
\theoremstyle{plain}

\newtheorem{lemma}{Lemma}[section]

\newtheorem{proposition}{Proposition}[section]

\numberwithin{equation}{section}

\begin{document}
\title[A Bernstein problem]{A Bernstein problem for special Lagrangian equations}
\author{Yu YUAN}
\address{Department of Mathematics\\
University of Chicago\\
5734 S. University Ave., Chicago, IL 60637\\
and University of Washington\\
Seattle, WA 98195}
\email{yuan@math.uchicago.edu}
\date{November 25, 2001. Submitted on December 10, 2001.}
\maketitle

\section{Introduction}

In this paper we derive a Bernstein type result for the special Lagrangian
equation%
\begin{equation}
F\left(  D^{2}u\right)  =\arctan\lambda_{1}+\cdots+\arctan\lambda
_{n}=c,\label{Eslg}%
\end{equation}
where $\lambda_{i}s$ are the eigenvalues of the Hessian $D^{2}u.$ Namely, any
global convex solution to (\ref{Eslg}) in $R^{n}$ must be a quadratic
polynomial. Recall the classical result, any global convex solution in $R^{n}$
to the Laplace equation $\triangle u=\lambda_{1}+\cdots+\lambda_{n}=c$ or the
Monge-Amp\`{e}re equation $\log\det D^{2}u=\log\lambda_{1}+\cdots+\log
\lambda_{n}=c$ must be quadratic.

Equation (\ref{Eslg}) originates from special Lagrangian geometry [HL]. The
(Lagrangian) graph $(x,\triangledown u\left(  x\right)  )\subset\mathbb{R}%
^{n}\times\mathbb{R}^{n}$ is called special when the argument of the complex
number $\left(  1+\sqrt{-1}\lambda_{1}\right)  \cdots\left(  1+\sqrt
{-1}\lambda_{n}\right)  $ is constant $c$ or $u$ satisfies (\ref{Eslg}), and
it is special if and only if $(x,\triangledown u\left(  x\right)  )$ is a
minimal surface in $\mathbb{R}^{n}\times\mathbb{R}^{n}$ [HL, Theorem 2.3,
Proposition 2.17].

In terms of minimal surface, our result is the following

\begin{theorem}
Suppose $M=(x,\triangledown u)$ is a minimal surface in $\mathbb{R}^{n}%
\times\mathbb{R}^{n}$ and $u$ is a smooth convex function in $\mathbb{R}^{n}$.
Then $M$ is a plane.
\end{theorem}

In fact, we have stronger results.

\begin{theorem}
Suppose $M=(x,\triangledown u)$ is a minimal surface in $\mathbb{R}^{n}%
\times\mathbb{R}^{n}$ and $u$ is a smooth function in $\mathbb{R}^{n}$ whose
Hessian satisfies $D^{2}u\geq-\epsilon(n)I,$ where $\varepsilon\left(
n\right)  $ is a small dimensional constant. Then $M$ is a plane.
\end{theorem}

\begin{theorem}
Suppose $M=(x,\triangledown u)$ is a minimal surface in $\mathbb{R}^{3}\times$
$\mathbb{R}^{3}$ and $u$ is a smooth function in $\mathbb{R}^{3}$ whose
Hessian satisfies $D^{2}u\geq-CI$. Then $M$ is a plane.
\end{theorem}

The lower bound on the Hessian $D^{2}u$ is necessary for Theorem 1.3, as one
sees from the following example. Let $u$ be a harmonic function in
$\mathbb{R}^{2},$ say, $u=x_{1}^{3}-3x_{1}x_{2}^{2},$ then $(x,\triangledown
u\left(  x\right)  )$ is a minimal surface in $\mathbb{R}^{4},$ which is not a plane.

Borisenko [Bo] proved Theorem 1.1 under the additional assumption that $u$
grows linearly at $\infty$ and $\arctan\lambda_{1}+\cdots+\arctan\lambda
_{n}=k\pi.$ For $c=k\pi,$ the special Lagrangian equation (\ref{Eslg}) 
in $R^{3}$ also takes the form%
\begin{equation}
\triangle u=\det D^{2}u. \label{Elm}%
\end{equation}
It was proved in [BCGJ] that any strictly convex solution to (\ref{Elm}) in
$\mathbb{R}^{3}$ with quadratic growth at $\infty$ must be quadratic.

Fu [F] showed that any global minimal surface $(x,\triangledown u\left(
x\right)  )\subset\mathbb{R}^{2}\times$ $\mathbb{R}^{2}$ is either a plane or
the potential $u$ is harmonic. This result also follows from Theorem 1.3
easily. We may assume $c\geq0$ in the special Lagrangian equation
$\arctan\lambda_{1}+\arctan\lambda_{2}=c.$ Then either $c=0,$ that is
$\triangle u=0$, or $\left(  D^{2}u\right) > -\frac{1}{\tan c}I,$ which in
turn implies that $u$ is quadratic by Theorem 1.3.

The heuristic idea of the proof of Theorem 1.1 is to find a subharmonic
function $S$ in terms of the Hessian $D^{2}u$ such that $S$ achieves its
maximum at a finite point in $R^{n}.$ By the strong maximum principle, $S$ is
constant. Consequently, $D^{2}u$ is a constant matrix. The right function $S$
is the one associated to the volume form of $M\;$in $R^{2n},$ $\det\left(
I+D^{2}uD^{2}u\right) ,$ see Lemma \ref{Lcomp}. However the nonnegative
Hessian $D^{2}u$ has no upper bound. We make a (Lewy) rotation of the $\left(
x,y\right)  \in\mathbb{R}^{n}\times\mathbb{R}^{n}$ coordinate system to
$\bar{x}=\left(  x+y\right)  /\sqrt{2},$ \ $\bar{y}=\left(  -x+y\right)
/\sqrt{2}.$ The special Lagrangian property of $M$ is invariant, and $M\;$has
a new representation $\left(  \bar{x},\triangledown\bar{u}\left(  \bar
{x}\right)  \right)  $ with the potential function $\bar{u}$ satisfying
$-I\leq\left(  D^{2}\bar{u}\right)  =\left(  I+D^{2}u\right)  ^{-1}\left(
-I+D^{2}u\right)  \leq I.$ To make the whole idea work, we need the machinery
from geometric measure theory, see Section 3.

Note that the special Lagrangian feature of the minimal surface $M=\left(
\bar{x},\triangledown\bar{u}\left(  \bar{x}\right)  \right)  $ is essential in
finding a subharmonic function. The function $\det\left(  I+D^{2}\bar{u}%
D^{2}\bar{u}\right)  $ is subharmonic as long as $-I\leq\left(  D^{2}\bar
{u}\right)  \leq I,$ in which case $\det\left(  I+D^{2}\bar{u}D^{2}\bar
{u}\right)  \leq2^{n}.$ In contrast, for general minimal surface $M=\left(
x,f\left(  x\right)  \right)  \subset R^{n}\times R^{k}$ with high
co-dimension $k\geq2,$ assuming that
\[
\det\left[  I+\left(  \triangledown f\right)  ^{t}\left(  \triangledown
f\right)  \right]  \leq K<\left[  \cos\left(  \pi/\left(  2\sqrt{2p}\right)
\right)  \right]  ^{-2p}%
\]
with $p=\min\left\{  n,k\right\}  ,$ Fischer-Colbrie [F-C] and Hildebrandt,
Jost, and Widman [HJW] were able to show that the composition of the square of
the distance function on the Grassmanian manifold $G\left(  n,k\right)  $ with
the harmonic map from $M\;$to $G\left(  n,k\right)  $ is subharmonic. Later
on, Jost and Xin [JX] proved the same thing under the assumption that
$\det\left[  I+\left(  \triangledown f\right)  ^{t}\left(  \triangledown
f\right)  \right]  \leq K<4.$ As a consequence, Bernstein type results were
obtained in all these papers.

Theorem 1.2 is just a consequence of Allard's $\varepsilon$-regularity theory,
once Theorem 1.1 is available.

Theorem 1.3 relies on the well-known result that any non-parametric minimal
cone of dimension three must be flat, see [F-C] and [B]. A quick ``PDE'' proof
of this fact was found in a recent paper [HNY]. Whether Theorem 1.3 holds true
in higher dimensional case remains an issue to us.

\textbf{Notation.} $\partial_{i}=\frac{\partial}{\partial x_{i}}%
,\;\partial_{ij}=\frac{\partial^{2}}{\partial x_{i}\partial x_{j}}%
,\;u_{i}=\partial_{i}u,\;u_{ji}=\partial_{ij}u,$ etc.

\section{Preliminary computations}

Let $\left(  x,\triangledown u\left(  x\right)  \right)  \subset\mathbb{R}%
^{n}\times\mathbb{R}^{n}$ be a non-parametric minimal surface, then we have%
\begin{equation}
\triangle_{g}\left(  x,\triangledown u\left(  x\right)  \right)  =0,\label{me}%
\end{equation}
where $\triangle_{g}=\sum_{i,j=1}^{n}\frac{1}{\sqrt{\det g}}\partial
_{i}\left(  \sqrt{\det g}g^{ij}\partial_{j}\right)  $ is the Laplace-Beltrami
operator of the induced metric $g=\left(  g_{ij}\right)  =\left(
I+D^{2}uD^{2}u\right)  $ with $\left(  g^{ij}\right)  =\left(  g_{ij}\right)
^{-1}.$ Notice that $\triangle_{g}x=0,$ $\triangle_{g}$ also takes the form%
\begin{equation}
\triangle_{g}=\sum_{i,j=1}^{n}g^{ij}\partial_{ij}.\label{ND}%
\end{equation}

\begin{lemma}
\label{Lcomp}Let $\left(  x,\triangledown u\left(  x\right)  \right)
\subset\mathbb{R}^{n}\times\mathbb{R}^{n}$ be a minimal surface. Suppose the
Hessian $D^{2}u$ of the smooth function $u$ is diagonalized at $p,$
$D^{2}u\left(  p\right)  =diag[\lambda_{1},\cdots,\lambda_{n}].$ Then%
\begin{equation}
\triangle_{g}\log\det g=\sum_{i,j=1}^{n}g^{ij}\partial_{ij}\log\det
g\overset{p}{=}\sum_{a,b,c=1}^{n}2g^{aa}g^{bb}g^{cc}u_{abc}^{2}\left(
1+\lambda_{b}\lambda_{c}\right)  .
\end{equation}
\end{lemma}

\begin{proof}
As preparation, we first compute the first and second order derivatives of the
metric $g.$%
\begin{align}
\partial_{j}g_{ab} &  =\sum_{k=1}^{n}\left(  u_{akj}u_{kb}+u_{ak}%
u_{kbj}\right)  \label{1dg}\\
&  \overset{p}{=}u_{abj}\left(  \lambda_{a}+\lambda_{b}\right)  .\nonumber
\end{align}%
\begin{align}
\partial_{i}g^{ab} &  =\sum_{k=1}^{n}-g^{ak}\partial_{i}g_{kl}g^{lb}%
\label{1dgi}\\
&  \overset{p}{=}-g^{aa}\partial_{i}g_{ab}g^{bb}\nonumber\\
&  \overset{p}{=}-g^{aa}g^{bb}u_{abi}\left(  \lambda_{a}+\lambda_{b}\right)
.\nonumber
\end{align}%

\begin{align*}
\partial_{ij}g_{ab} &  =\sum_{k=1}^{n}\left(  u_{akji}u_{kb}+u_{akj}%
u_{kbi}+u_{aki}u_{kbj}+u_{ak}u_{kbji}\right)  \\
&  \overset{p}{=}u_{abji}\left(  \lambda_{a}+\lambda_{b}\right)  +\sum
_{k=1}^{n}\left(  u_{akj}u_{kbi}+u_{aki}u_{kbj}\right)  .
\end{align*}
We need to substitute the 4$^{th}$ order derivative of $u$ with lower order
derivatives, we use the minimal surface equation (\ref{me}) with (\ref{ND}),%
\[
\triangle_{g}u_{a}=\sum_{i,j=1}^{n}g^{ij}\partial_{ij}u_{a}=0.
\]
Take the derivative with respect to $x_{b},$ we have%
\[
\sum_{i,j=1}^{n}\left(  g^{ij}\partial_{ij}u_{ab}+\partial_{b}g^{ij}%
\partial_{ij}u_{a}\right)  =0.
\]
Then%
\[
\sum_{i,j=1}^{n}g^{ij}\partial_{ij}u_{ab}\overset{p}{=}\sum_{i,j=1}^{n}%
g^{ii}g^{jj}u_{ijb}\left(  \lambda_{i}+\lambda_{j}\right)  u_{aji}%
\]
and%
\begin{equation}
\sum_{i,j=1}^{n}g^{ij}\partial_{ij}g_{ab}\overset{p}{=}\sum_{i,j=1}^{n}%
g^{ii}g^{jj}u_{ijb}u_{aji}\left(  \lambda_{i}+\lambda_{j}\right)  \left(
\lambda_{a}+\lambda_{b}\right)  +\sum_{i,k=1}^{n}2g^{ii}u_{aki}u_{kbi}%
.\label{2dg}%
\end{equation}
Relying on (\ref{1dg}) (\ref{1dgi}) (\ref{2dg}), we arrive at%
\begin{align*}
&  \sum_{i,j=1}^{n}g^{ij}\partial_{ij}\log\det g\\
&  =\sum_{i,j,a,b=1}^{n}g^{ij}\partial_{i}\left(  g^{ab}\partial_{j}%
g_{ab}\right)  \\
&  =\sum_{i,j,a,b=1}^{n}\left(  g^{ij}\partial_{i}g^{ab}\partial_{j}%
g_{ab}+g^{ij}g^{ab}\partial_{ij}g_{ab}\right)  \\
&  \overset{p}{=}\sum_{i,a,b=1}^{n}-g^{ii}g^{aa}g^{bb}u_{abi}^{2}\left(
\lambda_{a}+\lambda_{b}\right)  ^{2}+\sum_{i,j,a=1}^{n}2g^{aa}g^{ii}%
g^{jj}u_{aji}^{2}\left(  \lambda_{i}+\lambda_{j}\right)  \lambda_{a}%
+\sum_{i,k,a=1}^{n}2g^{aa}g^{ii}u_{aki}^{2}\\
&  \overset{p}{=}\sum_{a,b,c=1}^{n}-g^{aa}g^{bb}g^{cc}u_{abc}^{2}\left(
\lambda_{a}+\lambda_{b}\right)  ^{2}+\sum_{a,b,c=1}^{n}2g^{aa}g^{bb}%
g^{cc}u_{abc}^{2}\left(  \lambda_{b}+\lambda_{c}\right)  \lambda_{a}%
+\sum_{a,b,c=1}^{n}2g^{aa}g^{cc}u_{abc}^{2}\\
&  \overset{p}{=}\sum_{a,b,c=1}^{n}-2g^{aa}g^{bb}g^{cc}u_{abc}^{2}\left(
\lambda_{b}^{2}+\lambda_{a}\lambda_{b}\right)  +\sum_{a,b,c=1}^{n}%
4g^{aa}g^{bb}g^{cc}u_{abc}^{2}\lambda_{a}\lambda_{b}+\sum_{a,b,c=1}^{n}%
2g^{aa}g^{cc}u_{abc}^{2}\\
&  \overset{p}{=}\sum_{a,b,c=1}^{n}2g^{aa}g^{cc}u_{abc}^{2}\left(
-g^{bb}\lambda_{b}^{2}+1\right)  +\sum_{a,b,c=1}^{n}2g^{aa}g^{bb}g^{cc}%
u_{abc}^{2}\lambda_{a}\lambda_{b}\\
&  \overset{p}{=}\sum_{a,b,c=1}^{n}2g^{aa}g^{bb}g^{cc}u_{abc}^{2}\left(
1+\lambda_{a}\lambda_{b}\right)  \\
&  \overset{p}{=}\sum_{a,b,c=1}^{n}2g^{aa}g^{bb}g^{cc}u_{abc}^{2}\left(
1+\lambda_{b}\lambda_{c}\right)  ,
\end{align*}
where we use $g^{bb}\overset{p}{=}\frac{1}{1+\lambda_{b}^{2}}.$ This finishes
the proof of Lemma 2.1.
\end{proof}

\begin{proposition}
\label{Pslgc}Let $C=\left(  x,\triangledown u\left(  x\right)  \right)
\subset\mathbb{R}^{2n}$ be a minimal cone, smooth away from the origin.
Suppose the Hessian $D^{2}u$ satisfies $-I\leq\left(  D^{2}u\right)  \leq I.$
Then $C$ is a plane.
\end{proposition}

\begin{proof}
Since $\left(  x,\triangledown u\left(  x\right)  \right)  $ is cone,
$\triangledown u\left(  x\right)  $ is homogeneous degree one and
$D^{2}u\left(  x\right)  $ is homogeneous degree zero. It follows that
$\log\det g=\log\det\left(  I+D^{2}uD^{2}u\right)  $ takes its maximum at a
finite point (away from $0$) in $\mathbb{R}^{n}.$ By the assumption
$-I\leq\left(  D^{2}u\right)  \leq I,$ it follows from Lemma \ref{Lcomp} that%
\[
\sum_{i,j=1}^{n}g^{ij}\partial_{ij}\log\det g\geq0.
\]
By the strong maximum principle, we see that $\log\det g\equiv const..$
Applying Lemma \ref{Lcomp} again, we obtain%
\[
0\overset{p}{=}\sum_{a,b,c=1}^{n}2g^{aa}g^{bb}g^{cc}u_{abc}^{2}\left(
1+\lambda_{b}\lambda_{c}\right)  \geq0.
\]
Then%
\[
u_{abc}^{2}\left(  1+\lambda_{a}\lambda_{b}\right)  =u_{abc}^{2}\left(
1+\lambda_{b}\lambda_{c}\right)  =u_{abc}^{2}\left(  1+\lambda_{c}\lambda
_{a}\right)  =0.
\]
Observe that one of $\lambda_{a}\lambda_{b},\;\ \lambda_{b}\lambda_{c},$
and\ $\lambda_{c}\lambda_{a}$ must be nonnegative, we get $u_{abc}\left(
p\right)  =0.$ Since the point $p$ in Lemma \ref{Lcomp} can be arbitrary, we
conclude that $D^{3}u\equiv0.$ Consequently, $u$ is a quadratic function and
the cone $\left(  x,\triangledown u\left(  x\right)  \right)  $ is a plane.
\end{proof}

\section{Proof of theorems}

\begin{proof}
[Proof of Theorem 1.1]Step A. We first seek a better representation of $M$ via
Lewy transformation. We rotate the $\left(  x,y\right)  \in\mathbb{R}%
^{n}\times\mathbb{R}^{n}$ coordinate system to $\left(  \bar{x},\bar
{y}\right)  $ by $\pi/4,$ namely, set $\bar{x}=\left(  x+y\right)  /\sqrt{2},$
\ $\bar{y}=\left(  -x+y\right)  /\sqrt{2}.$ Then $M$ has a new parametrization%
\[
\left\{
\begin{array}
[c]{c}%
\bar{x}=\frac{1}{\sqrt{2}}\left(  x+\triangledown u\left(  x\right)  \right)
\\
\bar{y}=\frac{1}{\sqrt{2}}\left(  -x+\triangledown u\left(  x\right)
\right)
\end{array}
\right.  .
\]
Since $u$ is convex, we have%
\begin{align*}
\left|  \bar{x}^{2}-\bar{x}^{1}\right|^{2}   &  =\frac{1}{2}\left[  \left|
x^{2}-x^{1}\right|  ^{2}+2\left(  x^{2}-x^{1}\right)  \cdot\left(
\triangledown u\left(  x^{2}\right)  \ -\triangledown u\left(  x^{1}\right)
\ \right)  +\left|  \triangledown u\left(  x^{2}\right)  \ -\triangledown
u\left(  x^{1}\right)  \ \right|  ^{2}\right]  \\
&  \geq\frac{1}{2}\left|  x^{2}-x^{1}\right|  ^{2}.
\end{align*}
It follows that $M$ is still a graph over the whole $\bar{x}$ space
$\mathbb{R}^{n}.$ Further $M$ is still a Lagrangian graph over $\bar{x},$ that
means $M$ has the representation $\left(  \bar{x},\triangledown\bar{u}\left(
\bar{x}\right)  \right)  $ with a potential function $\bar{u}\in C^{\infty
}\left(  \mathbb{R}^{n}\right)  $ (cf. [HL, Lemma 2.2]).

Note that any tangent vector to \ $M$ takes the form%
\[
\frac{1}{\sqrt{2}}\left(  \left(  I+D^{2}u\left(  x\right)  \right)  e,\left(
-I+D^{2}u\left(  x\right)  \right)  e\right)  ,
\]
where $e\in\mathbb{R}^{n}.$ It follows that%
\[
D^{2}\bar{u}\left(  \bar{x}\right)  =\left(  I+D^{2}u\left(  x\right)
\right)  ^{-1}\left(  -I+D^{2}u\left(  x\right)  \right)  .
\]
By the convexity of $u,$ we have%
\[
-I\leq\left(  D^{2}\bar{u}\right)  \leq I.
\]

Step B. The remaining proof is routine. We ``blow down'' $M\;$at $\infty.$
Without loss of generality, we assume $\bar{u}\left(  0\right)
=0,\;\triangledown\bar{u}\left(  0\right)  =0.$ Set $M_{k}=\left(  \bar
{x},\triangledown\bar{u}_{k}\right)  ,$ where%
\[
\bar{u}_{k}\left(  \bar{x}\right)  =\frac{\bar{u}\left(  k\bar{x}\right)
}{k^{2}},\;\;\;k=1,2,3,\ \cdots.
\]
We see that $M_{k}$ is still a minimal surface and $-I\leq\left(
D^{2}\bar{u}_{k}\right)  \leq I.\;$Then there exists a subsequence, still
denoted by $\left\{  \bar{u}_{k}\right\}  $ and $v\in C^{1,1}\left(
R^{n}\right)  $ such that%
\[
\bar{u}_{k}\rightarrow v\;\;\text{in \ \ }C_{loc}^{1,\alpha}\left(
R^{n}\right)
\]
and%
\[
-I\leq\left(  D^{2}v\right)  \leq I.
\]
We apply the compactness theorem (cf. [S, Theorem 34.5] to conclude that
$M_{v}=\left(  \bar{x},\triangledown v\left(  \bar{x}\right)  \right)  $ is a
minimal surface, By the monotonicity formula (cf. [S, p.84]) and Theorem 19.3
in [S], we know that $M_{v}$ is a minimal cone.

We claim that $M_{v}$ is smooth away from the vertex. Suppose $M_{v}$ is
singular at $P$ away from the vertex. We blow up $M_{v}$ at $P$ to get a
tangent cone, which is a lower dimensional special Lagrangian cone cross a
line, repeat the procedure if the resulting cone is still singular away from
the vertex. Finally we get a special Lagrangian cone which is smooth away from
the vertex, and the eigenvalues of the Hessian of the potential function are
bounded between $-1$ and $1.$ By Proposition \ref{Pslgc}, the cone is flat.
This is a contradiction to Allard's regularity result (cf. [S, Theorem 24.2]).

Applying Proposition \ref{Pslgc} to $M_{v},$ we see that $M_{v}$ is flat.

Step C. By our blow-down procedure and the monotonicity formula, we see that%
\[
\lim_{r\rightarrow+\infty}\frac{\mu\left(  \frak{B}_{r}\left(  0,0\right)
\cap M\right)  }{\left|  B_{r}\right|  }=1,
\]
where $B_{r}$ is the ball with radius $r$ in $\mathbb{R}^{n},$ $\frak{B}%
_{r}\left(  0,0\right)  $ is the ball with radius $r$ and center $\left(
0,0\right)  $ in $\mathbb{R}^{n}\times\mathbb{R}^{n},$ and $\mu\left(
\frak{B}_{r}\left(  0,0\right)  \cap M\right)  $ is the area of $M$ inside
$\frak{B}_{r}\left(  0,0\right)  .$ Since $M$ is smooth, we have%
\[
\lim_{r\rightarrow0}\frac{\mu\left(  \frak{B}_{r}\left(  0,0\right)  \cap
M\right)  }{\left|  B_{r}\right|  }=1.
\]
Consequently, for $r_{2}>r_{1}>0,$\ the monotonicity formula reads
\[
0=\frac{\mu\left(  \frak{B}_{r_{2}}\left(  0,0\right)  \cap M\right)
}{\left|  B_{r_{2}}\right|  }-\frac{\mu\left(  \frak{B}_{r_{1}}\left(
0,0\right)  \cap M\right)  }{\left|  B_{r_{1}}\right|  }=\int_{\frak{B}%
_{r_{2}}\backslash\frak{B}_{r_{1}}}\frac{\left|  D^{\perp}r\right|  ^{2}%
}{r^{n}}d\mu,
\]
where $r=\left|  \left(  x,y\right)  \right|  ,$ $D^{\perp}r$ is the
orthogonal projection of $Dr$ to the normal space of $M,$ and $d\mu$ is the
area form on $M.$ Therefore, we see that $M\;$is a plane.
\end{proof}

\textbf{Remark}. In Step B, we use the heavy compactness result (cf. [S,
Theorem 34.5]) just for a short presentation of the proof. One can also take
advantage of the special Lagrangian equation (\ref{Eslg}), use the compactness
result for viscosity solution to derive that $M_{v}=\left(  \bar
{x},\triangledown v\left(  \bar{x}\right)  \right)  $ is a minimal surface,
see Lemma 2.2 in [Y].

\begin{proposition}
\label{Pslgp}There exist a dimensional constant $\varepsilon^{\prime}\left(
n\right)  >0$ such that any minimal surface $\left(  x,\triangledown u\left(
x\right)  \right)  \subset\mathbb{R}^{n}\times\mathbb{R}^{n}$ with $-\left(
1+\varepsilon^{\prime}\left(  n\right)  \right)  I\leq\left(  D^{2}u\right)
\leq\left(  1+\varepsilon^{\prime}\left(  n\right)  \right)  I$ for
$x\in\mathbb{R}^{n}$, must be a plane.
\end{proposition}

\begin{proof}
Suppose not. Then there exists a sequence of minimal surface $M_{k}=\left(
x,\triangledown u_{k}\right)  \subset\mathbb{R}^{n}\times\mathbb{R}^{n}$ such
that $-\left(  1+\frac{1}{k}\right)  I\leq\left(  D^{2}u_{k}\right)
\leq\left(  1+\frac{1}{k}\right)  I$ and $M_{k}$ is not a plane. By Allard's
regularity result (cf. [S, Theorem 24.2]) the density $D_{k}$ for $M_{k}$
satisfies%
\[
D_{k}\geq1+\delta\left(  n\right)  ,
\]
where $\delta\left(  n\right)  >0$ is a dimensional constant and%
\[
D_{k}=\lim_{r\rightarrow+\infty}\frac{\mu\left(  \frak{B}_{r}\cap
M_{k}\right)  }{\left|  B_{r}\right|  }.
\]
By a similar argument as Step B in the proof of Theorem 1.1, we extract a
subsequence of $\left\{  v_{k}\right\}  $ converging to $V_{\infty}$ in
$C_{loc}^{1,\alpha}\left(  R^{n}\right)  $ such that $M_{\infty}=\left(
x,\triangledown V_{\infty}\left(  x\right)  \right)  $ is a smooth minimal
surface in $\mathbb{R}^{n}\times\mathbb{R}^{n}$ with $-I\leq\left(
D^{2}u_{\infty}\right)  \leq I$ and $D_{\infty}\geq1+\delta\left(  n\right)
.$ By our Theorem 1.1, $M_{\infty}$ is a plane and $D_{\infty}=1.$ This
contradiction finishes the proof of the proposition.
\end{proof}

\begin{proof}
[Proof of Theorem 1.2]We repeat the rotation argument in Step A of the proof
of Theorem 1.1 to get a new representation for $M,$ $\left(  \bar
{x},\triangledown\bar{u}\left(  \bar{x}\right)  \right)  $ with%
\[
-\left(  1+\frac{2\varepsilon\left(  n\right)  }{1-\varepsilon\left(
n\right)  }\right)  I\leq\left(  D^{2}\bar{u}\right)  \leq I.
\]
We choose $\varepsilon\left(  n\right)  =\frac{\varepsilon^{\prime}\left(
n\right)  }{2+\varepsilon^{\prime}\left(  n\right)  }$ and apply Proposition
\ref{Pslgp}. Then Theorem 1.2 follows.
\end{proof}

\begin{proof}
[Proof of Theorem 1.3]The strategy is similar to the proof of Theorem 1.1.

Step A. We first make a different rotation of the coordinate system to get a
better representation of $M.$ Set $\bar{x}=\frac{1}{\sqrt{1+4C^{2}}}\left(
2Cx+y\right)  ,$ $\;\bar{y}=\frac{1}{\sqrt{1+4C^{2}}}\left(  -x+2Cy\right)  .$
Then $M$ has a new parametrization%
\[
\left\{
\begin{array}
[c]{c}%
\bar{x}=\frac{1}{\sqrt{1+4C^{2}}}\left(  2Cx+\triangledown u\left(  x\right)
\right)  \\
\bar{y}=\frac{1}{\sqrt{1+4C^{2}}}\left(  -x+2C\triangledown u\left(  x\right)
\right)
\end{array}
\right.  .
\]
Since $u+\frac{1}{2}C\left|  x\right|  ^{2}$ is convex, we have%
\begin{align*}
\left|  \bar{x}^{2}-\bar{x}^{1}\right|^{2}   &  =\frac{1}{1+4C^{2}}\left[
\begin{array}
[c]{c}%
C^{2}\left|  x^{2}-x^{1}\right|  ^{2}+2C\left(  x^{2}-x^{1}\right)
\cdot\left(  \triangledown u\left(  x^{2}\right)  +Cx^{2}-\triangledown
u\left(  x^{1}\right)  -Cx^{1}\right)  \\
+\left|  \triangledown u\left(  x^{2}\right)  +Cx^{2}-\triangledown u\left(
x^{1}\right)  -Cx^{1}\right|  ^{2}%
\end{array}
\right]  \\
&  \geq\frac{1}{1+4C^{2}}C^{2}\left|  x^{2}-x^{1}\right|  ^{2}.
\end{align*}
As in the proof of Theorem 1.1, we get a new representation for $M=\left(
\bar{x},\triangledown\bar{u}\left(  \bar{x}\right)  \right)  $ and%
\[
D^{2}\bar{u}\left(  \bar{x}\right)  =\left(  2CI+D^{2}u\right)  ^{-1}\left(
-I+2CD^{2}u\left(  x\right)  \right)  .
\]
From $D^{2}u\geq-CI,$ we see that%
\[
-\frac{1+2C^{2}}{C}I\leq\left(  D^{2}\bar{u}\right)  \leq2CI.
\]

Step B. As step B in the proof of Theorem 1.1, any tangent cone of $M$ at
$\infty$ is flat. The only difference is that, instead of relying on
Proposition \ref{Pslgc}, we use the fact that any non-parametric minimal cone
of dimension three must be flat, see [F-C, Theorem 2.3], [B, Theorem]. For a
quick PDE proof of this fact, see [HNY, p.2].

Step C is exactly as in the proof of Theorem 1.1.

Therefore, we conclude Theorem 1.3.
\end{proof}

\begin{center}
\textbf{Acknowledgements}
\end{center}

The author thanks Dan Freed and Karen Uhlenbeck for ``forcing'' him to give a
talk on minimal surfaces in the GADET seminar at University of Texas at
Austin, which resulted the present work. The author is grateful to Luis
Caffarelli for discussions and pointing out Theorem 1.2, to Gang Tian for
pointing out the simple argument in Step C of the proof of Theorem 1.1, which
avoids the usual approach via Allard's regularity. The author is supported in
part by an NSF grant.

\end{document}